\documentclass{amsart}

\newtheorem{theorem}{Theorem}[section]
\newtheorem{lemma}[theorem]{Lemma}
\newtheorem{proposition}[theorem]{Proposition}

\theoremstyle{definition}
\newtheorem{definition}[theorem]{Definition}
\newtheorem{example}{Example}[section]

\theoremstyle{remark}
\newtheorem{remark}{Remark}[section]

\numberwithin{equation}{section}

\newcommand{\abs}[1]{{\left|#1\right|}}

\newcommand{\norm}[2][p]{{\parallel#2\parallel}_{#1}}

\newcommand{\im}{\mathop{\mathrm {Im}}\nolimits}
\newcommand{\re}{\mathop{\mathrm {Re}}\nolimits}
\newcommand{\C}{{\mathbb C}}
\newcommand{\R}{{\mathbb R}}

\newcommand{\D}{{\mathcal D}}
\newcommand{\E}{{\mathcal E}}
\newcommand{\Rect}{{\mathcal R}}
\newcommand{\K}{{\mathcal K}}

\begin{document}

\title[Weighted spaces in the upper half plane]{Weighted Banach spaces of
holomorphic functions in the upper half plane}

\author{Martin At. Stanev}

\curraddr{Department of Mathematics,
University  of Chemical Techonology  and Metallurgy,
8, bulv. Kliment Ohridski,
 Sofia 1756, Bulgaria}
\email{stanevm@adm1.uctm.edu  or  stanevm@Math1.uctm.edu }
\thanks{}


\date{}


\begin{abstract}
We define a new class of weighted Banach spaces of holomorphic functions in
the upper half plane. Three basical theorems on these spaces are proved. Two
of them investigate necessary and sufficient conditions on the weight
function in  order to obtain that the defined spaces are not the trivial
ones. The third theorem deals with the behaviour at infinity of some
holomorphic functions. These three theorems all are new and do not represent
neither generalizations nor particular cases of  any previous results.
\end{abstract}

\maketitle

\section{Introduction}\label{First}

Let $\C$ be the complex plane and  let $\R$ be  the real axis. We use the
following  notations: the real part of $z\in\C$ is denoted by $\re z$; the
imaginary  part of $z\in\C$ is denoted by $\im  z$; and $\D$ stands for the
upper half  plane of $\C$, i.e. $\D=\{z|\,z\in\C,\im z>0\}$.
\begin{definition}\label{Def}
A function $p:(0,\infty)\to(0,\infty)$ is said to be  weight function  if
and  only if
\begin{equation}\label{Positive}
\inf_{t,t\in[\frac1c,c]}p(t)>0,\qquad \forall c>1.
\end{equation}
Let $p$ be any weight function. We introduce the notation
\[
\norm{f}=\sup_{z,z\in \D}p(\im z)\abs{f(z)},
\]
where $f$ is any holomorphic function in the upper half plane $\D$.
 And,  we define the spaces $\Lambda(p)$  and $\lambda(p)$ as follows:
\begin{itemize}
\item $\Lambda(p)$ is a normed complex function space and $f\in \Lambda(p)$
if and only if $f$ is such a holomorphic function in the upper half plane
$\D$ that $\norm{f}<\infty$, and $\norm{f}$ is the norm of $f\in\Lambda(p)$;
\item $\lambda(p)$ is a normed complex function space and $f\in\lambda(p)$
if and only if $f$ is a holomorphic function in $\D$ such that satisfies the
following two conditions:
\begin{enumerate}
   \item $\norm{f}<\infty$,
   \item  for every positive number $\varepsilon,\varepsilon>0$, there
exists a number $c,c\ge1$ such that
\[
\sup_{z,z\in \D\setminus {\K}_c}p(\im z)\abs{f(z)}\le\varepsilon,
\]
where ${\K}_c=\{z|\,z\in \D, \frac1c\le \im z\le c,\abs{\re z}\le c\},c\ge
1$,
\end{enumerate}
and $\norm{f}$ is the norm of $f\in\lambda(p)$.
\end{itemize}
The spaces $\Lambda(p)$  and  $\lambda(p)$ are called  weighted Banach
spaces of holomorphic functions in the upper half plane.
\end{definition}

Here, in the present  Section, the main  definition and the  main results
are stated.  Theorems~\ref{Big} and~\ref{Small} solve completely  the
problem whether  the weighted Banach spaces  are trivial or not   by giving
the corresponding necessary  and sufficient conditions on  the weight
function.  Theorem~\ref{SmallBehaviour} concerns the  behaviour   of some
functions of~$\lambda(p)$ at infinity.

In Section~\ref{PR} there are examples of such spaces  and it is proved that
the weighted Banach spaces of holomorphic functions in the  upper half plane
all are Banach spaces.

In  Section~\ref{Proofs}  we prove our main results.

In this paper the main results are the following  three theorems.

\begin{theorem}\label{Big}
Let the function $p:(0,\infty)\to(0,\infty)$ meet the
condition~(\ref{Positive}).
Then the space $\Lambda(p)\ne\{0\}$ if and only if there are two real
numbers $a,b$ such that
\[
(-1)\ln p(t)\ge at+b,\qquad \forall t>0.
\]
\end{theorem}

\begin{theorem}\label{Small}
Let the function $p:(0,\infty)\to(0,\infty)$ meet the
condition~(\ref{Positive}).
Then the space $\lambda(p)\ne\{0\}$ if and  only if the following two
conditions on the weight function $p$ are fulfilled:
\begin{enumerate}
\item there are two real numbers $a,b$ such  that
\[
(-1)\ln p(t)\ge at+b,\qquad \forall t>0;
\]
\item $\displaystyle\lim_{t\to 0{+}}p(t)=0$.
\end{enumerate}
\end{theorem}

Further, for  brevity,   by $Mf(y)$ we denote
\[
Mf(y)=\sup_{z,\im z=y}\abs{f(z)},\qquad \forall y>0,f\in \Lambda(p),
\]
where $\Lambda(p)$ is any weighted  Banach  space of holomorphic functions
in the upper half plane.

\begin{theorem}\label{SmallBehaviour}
Let the function $p:(0,\infty)\to(0,\infty)$ be such that
\begin{eqnarray}
\inf_{t,t\in[\frac1c,c]}p(t)>0,&\qquad \forall c>1,\nonumber\\
\lambda(p)\ne\{0\}.&\nonumber
\end{eqnarray}
If  the holomorphic  function $f\in\lambda(p)\setminus\{0\}$ is such that
\[
\liminf_{t\to\infty}\frac{\ln Mf(t)}{t}<\infty,
\]
then for the limit value $a=\liminf\limits_{t\to \infty}\frac{\ln Mf(t)}{t}$
we obtain $a>-\infty$ and
\[
\lim_{t\to \infty}(\ln Mf(t)-at)=-\infty.
\]
\end{theorem}

\section{Preliminary results}\label{PR}

\subsection{Examples}

We consider some examples of weighted Banach spaces of holomorphic functions
in the upper half plane.
\begin{example}
Let $p(t)=1$, $\forall t>0$. Then the condition~(\ref{Positive}) is
fulfilled. Accordingly to the Definition~\ref{Def} the spaces $\Lambda(p)$
and $\lambda(p)$ are well-defined. In particular,
\[
\norm{f}=\sup_{z,\im z>0}\abs{f(z)},
\]
and $\Lambda(p)=H^\infty$. The assertion $\lambda(p)=\{0\}$ is a consequence
of Theorem~\ref{Small}.
\end{example}

\begin{example}
Let $p(t)=t$, $\forall t>0$. Then the condition~(\ref{Positive}) is
fulfilled. Accordingly to the Definition~\ref{Def} the spaces $\Lambda(p)$
and $\lambda(p)$ are well-defined. In particular,
\[
\norm{f}=\sup_{z,\im  z>0}\im z\abs{f(z)}.
\]
In order to charcterize the spaces $\Lambda(p),\lambda(p)$, we first recall
the definition of Bloch spaces $B$  and $B_0$ in the open unit disk (for
this  definition see~\cite{Anderson+Cl+Pomm1974})
\[
\begin{array}{rl}
f\in B\iff&\mbox{$f$ is holomorphic function}\\
&\mbox{in the open unit disk such that}\\
&\displaystyle \sup_{w,\abs{w}<1}(1-\abs{w}^2)\abs{f'(w)}<\infty;\\[6pt]
f\in B_0\iff&\displaystyle f\in B, \lim_{\abs{w}\to
1{-}}(1-\abs{w}^2)\abs{f'(w)}=0,\\
\end{array}
\]
where $f'$ stands for the derivative of  the function $f$, and
$\lim\limits_{\abs{w}\to  1{-}}$ is uniform with respect to the argument of
$w$.
We  set
\[
g_f(z)=4\frac{f'(\frac{1+iz}{1-iz})}{(1-iz)^2},\qquad f\in B,\im z>0,
\]
Then by a simple computation we obtain the following characterization of the
spaces $\Lambda(p)$ and $\lambda(p)$:
\begin{eqnarray}
f\in B\iff&g_f\in\Lambda(p),\nonumber\\
f\in B_0\iff&g_f\in\lambda(p).\nonumber
\end{eqnarray}
\end{example}

\begin{example}
Let $p(t)=t+e^t$,  $\forall t>0$. Then the condition~(\ref{Positive}) is
fulfilled. Accordingly to the Definition~\ref{Def} the spaces $\Lambda(p)$
and $\lambda(p)$ are well-defined. In particular,
\[
\norm{f}=\sup_{z,\im z>0}(\im z+e^{\im z})\abs{f(z)}.
\]
The space $\Lambda(p)$ does not contain any constant functions, and contains
the function $h$,
\[
h(z)=e^{2iz},\qquad \im z>0.
\]
Thus, $\Lambda(p)$ is a  proper subspace of  the  corresponding  spaces  of
previous two examples. The assertion $\lambda(p)=\{0\}$ is a consequence  of
Theorem~\ref{Small}.
\end{example}

\begin{example}
Let $p(t)=e^{t^2}$, $\forall t>0$.  Then the condition~(\ref{Positive}) is
fulfilled. Accordingly to the Definition~\ref{Def} the spaces $\Lambda(p)$
and $\lambda(p)$ are well-defined. In particular,
\[
\norm{f}=\sup_{z,\im z>0}e^{(\im z)^2}\abs{f(z)}.
\]
The assertions $\Lambda(p)=\{0\}$, $\lambda(p)=\{0\}$ are consequences
respectively of Theorems~\ref{Big} and~\ref{Small}.
\end{example}

\subsection{The completeness of $\Lambda(p)$ and $\lambda(p)$}

\begin{proposition}\label{Completeness}
Let the function $p:(0,\infty)\to(0,\infty)$ meet  the
condition~(\ref{Positive}). Then $\Lambda(p)$ and $\lambda(p)$,  both  are
Banach spaces.
\end{proposition}

\begin{remark}
The main idea of  the proof of Proposition~\ref{Completeness} is
demonstrated in~\cite{Anderson+Cl+Pomm1974}.
\end{remark}

\begin{proof}[Proof of proposition~\ref{Completeness}]
It is sufficient to prove proposition~\ref{Completeness} with additional
assumption that $\Lambda(p)\ne\{0\}$ and resp. $\lambda(p)\ne\{0\}$. Thus
let us assume the  corresponding one.

For brevity, we set ${\K}_c=\{z|\,z\in \D, \frac1c\le \im z\le c,\abs{\re
z}\le c\}$, $c\ge1$.

First we shell prove that $\Lambda(p)$ is complete.

Suppose $\{f_n\}_{n=1}^\infty$ be any  Cauchy's sequence in $\Lambda(p)$.
Then for every positive number $\varepsilon$ and for  every $c\ge 1$ there
is a number $n_0$ such that $\forall m,n>n_0$
\begin{equation}\label{eq}
\sup_{z,z\in {\K}_c}p(\im
z)\abs{f_n(z)-f_m(z)}\le\norm{f_n-f_m}\le\varepsilon.
\end{equation}
Hence, in particular,
\[
\sup_{z,z\in
{\K}_c}\abs{f_n(z)-f_m(z)}\le\frac{\varepsilon}{\inf\limits_{t\in
[\frac1c,c]}p(t)}.
\]
Thus we obtain that the  sequence $\{f_n\}_{n=1}^\infty$ is uniformly
convergent on any compact in  the  upper half plain $\D$. Then there is  a
holomorphic function $f$ in $\D$ which is  the uniform limit of
$\{f_n\}_{n=1}^\infty$ on every compact in $\D$. Therefore, by the
inequality~(\ref{eq}) it follows
\[
\sup_{z,z\in {\K}_c}p(\im z)\abs{f_n(z)-f(z)}\le\varepsilon.
\]
Hence, accordingly to  the choise of $c$ we obtain \[
\norm{f_n-f}\le\varepsilon,\qquad n>n_0.
\]
In particular, $f\in\Lambda(p)$ and
$\lim\limits_{n\to\infty}\norm{f_n-f}=0$.

Thus we obtain that $\Lambda(p)$ is complete and hence it is a Banach space.

Second, we shell prove that $\lambda(p)$ is complete.

Note $\lambda(p)$ is a linear subspace of the Banach space $\Lambda(p)$. So,
to obtain  the completeness of $\lambda(p)$ it is enough to prove that
$\lambda(p)$ is a closed subspace of $\Lambda(p)$.

We claim that $\lambda(p)$  is a closed subspace of $\Lambda(p)$.

Suppose $\{f_n\}_{n=1}^\infty$ is such  a sequence of elements of
$\lambda(p)$ that is convergent in the space $\Lambda(p)$. Let
$f\in\Lambda(p)$ be such that  $\lim\limits_{n\to \infty}\norm{f_n-f}=0$.

Then for every positive number $\varepsilon$ there is a number $n_0$ such
that
\[
\sup_{z,\im z>0}p(\im z)\abs{f_n(z)-f(z)}\le\frac{\varepsilon}2, \qquad
\forall n>n_0.
\]
Fix an $n$, $n>n_0$. Then it follows by $f_n\in\lambda(p)$ that there is
$c$, $c\ge1$, such that
\[
\sup_{z,z\in \D\setminus {\K}_c}p(\im z)\abs{f_n(z)}\le\frac{\varepsilon}2.
\]
Hence
\[
\sup_{z,z\in \D\setminus {\K}_c}p(\im z)\abs{f(z)}\le\varepsilon.
\]
and therefore $f\in\lambda(p)$. Thus the assertion which we claim above is
proved.
\end{proof}

\subsection{A preliminary lemma}

\begin{lemma}\label{Convexity}
Let the function $p:(0,\infty)\to(0,\infty)$ be such that the following  two
conditions are fulfilled:
\begin{eqnarray}
&\displaystyle \inf_{t,t\in[\frac1c,c]}p(t)>0,&\qquad\forall
c\ge1,\nonumber\\
&\Lambda(p)\ne\{0\}.&{}\nonumber
\end{eqnarray}
If $f\in\Lambda(p)\setminus\{0\}$ then
\begin{enumerate}
\item $\ln Mf$ is convex in $(0,\infty)$;
\item there  are two  numbers $a,b$ such that
\[
\ln Mf(t)\ge at+b, \qquad \forall t>0.
\]
\end{enumerate}
\end{lemma}

\begin{proof}[Proof of lemma~\ref{Convexity}]
Let $f$, $f\in\Lambda(p)\setminus\{0\}$, be an arbitrary non zero  element
of $\Lambda(p)$. Then
\[
\sup_{t,t>0}p(t)Mf(t)=\norm{f}<\infty.
\]
Hence
\[
\sup_{t,t\in[\frac1c,c]}Mf(t)\le\frac{\norm{f}}{\inf\limits_{t\in
[\frac1c,c]}p(t)}<\infty,\qquad \forall c\ge 1,
\]
and in particular the holomorphic in the upper half plane  $\D$ function $f$
is bounded in the band
\[
\{z|\,z\in \D,\frac1c\le \im z\le c\}, \qquad \forall  c\ge 1.
\]
Then by the Phragm\'en-Lindel\"of principle we obtain that  $\ln Mf$ is
convex in $[\frac1c,c]$, $\forall c\ge1$. Thus $\ln Mf$ is convex in
$(0,\infty)$.

So, the first assertion of Lemma~\ref{Convexity} is proved.
In order to prove the second assertion of this Lemma we set
\[
\E=\{(u,v)|\,u>0,v>\ln Mf(u)\}.
\]
From the  first assertion of Lemma~\ref{Convexity} it follows that  $\E$ is
an open convex subset of Euclidean plane $\{(u,v)|u\in\R,v\in\R\}$.
Further, by $f\ne 0$ it follows that there is $t$, $t>0$, such that $\ln
Mf(t)\ne{-}\infty$ (in  fact this inequality holds for every $t>0$, but  we
do not use it). Hence there exists a point $(u_0,v_0)$, $u_0>0$, $v_0\in\R$,
which does  not belong  to $\E$.  Then there  exists a line that
distinguishes the set $\E$ and the point $(u_0,v_0)$.

Hence, there are two numbers $a,b$ such that
\[
\ln Mf(t)\ge  at+b,\qquad\forall t>0.
\]
So, Lemma~\ref{Convexity} is proved.
\end{proof}

\section{Proof of main results}\label{Proofs}

\begin{proof}[Proof  of Theorem~\ref{Big}]
First we prove that if $\Lambda(p)\ne\{0\}$ then there are two numbers $a,b$
such that
\[
(-1)\ln p(t)\ge at+b, \qquad \forall t>0.
\]
Suppose $\Lambda(p)\ne\{0\}$ and let $f$, $f\in\Lambda(p)\setminus\{0\}$ be
an arbitrary non zero  function that belongs to $\Lambda(p)$. Then by
Lema~\ref{Convexity}, there exist two numbers $a_0$, $b_0$ such  that
\[
\ln Mf(t)\ge a_0t+b_0,\qquad\forall t>0.
\]
Further, by $f\in\Lambda(p)$ and  $\displaystyle
p(t)Mf(t)\le\norm{f}$, $\forall t>0$, it follows
\[
\begin{array}{rll}
(-1)\ln p(t)&\ge(-1)\ln\norm{f} +\ln Mf(t)&\\
&\ge(-1)\ln\norm{f}+a_0t+b_0,&\forall t>0.
\end{array}
\]
We set $a=a_0$, $b=b_0-\ln\norm{f}$ and obtain
\[
(-1)\ln p(t)\ge at+b,\qquad\forall t>0.
\]
Thus we prove that if $\Lambda(p)\ne\{0\}$ then the condition on  the
function $p$ which is stated in Theorem~\ref{Big} is fulfilled.

With view to show  $\Lambda(p)\ne\{0\}$, suppose that the function $p$
satisfies condition~(\ref{Positive}) and in addition there exist two numbers
$a$ and  $b$ such that
\[
(-1)\ln p(t)\ge at+b,\qquad\forall t>0.
\]
Then for such numbers  $a$, $b$ we set $f(z)=e^{-iaz+b}$, where $z\in\D$.
Then a direct  computation shows $\norm{f}\le 1$. So,
$f\in\Lambda(p)\setminus\{0\}$.

Thus the proof of Theorem~\ref{Big} is completed.
\end{proof}

\begin{proof}[Proof of Theorem~\ref{Small}]
Suppose $\lambda(p)\ne\{0\}$ and let $f$, $f\in\lambda(p)\setminus\{0\}$ be
an arbitrary non zero  function that belongs to $\lambda(p)$. In particular,
from  the inclusion $\lambda(p)\subset\Lambda(p)$ it follows
$f\in\Lambda(p)\setminus\{0\}$. Therefore, by Theorem~\ref{Big} the  first
condition on the function $p$ is fulfilled. Furthermore, according to
Definition~\ref{Def}
\[
\lim_{t\to 0{+}}p(t)Mf(t)=0.
\]
Moreover, by  Lemma~\ref{Convexity}, there exist two numbers $a_0$, $b_0$
such  that
\[
\ln Mf(t)\ge a_0t+b_0,\qquad\forall t>0.
\]
Hence,
\[
\begin{array}{rl}
0&\le\displaystyle\liminf_{t\to 0{+}}p(t)\le\limsup_{t\to 0{+}}p(t)\\[6pt]
&\le\displaystyle\limsup_{t\to 0{+}}p(t)e^{at+b}e^{-(at+b)}\\[6pt]
&\le\displaystyle\limsup_{t\to 0{+}}p(t)Mf(t)\lim_{t\to 0{+}}e^{-(at+b)}=0
\end{array}
\]
which yields $\lim\limits_{t\to 0{+}}p(t)=0$.

Thus, as we assert in Theorem~\ref{Small}, both conditions on function $p$
are  fulfilled.

Further, suppose that the function $p$ meets the condition~(\ref{Positive})
and in addition the following two conditions are  fulfilled, too:
\begin{enumerate}
\item there exist two numbers  $a$, $b$ such that
\[
(-1)\ln p(t)\ge at+b,\qquad\forall t>0.
\]
\item $\displaystyle\lim_{t\to0{+}}p(t)=0$.
\end{enumerate}
We claim that $\lambda(p)\ne\{0\}$.

Indded, the  function $f$ defined by
\[
f(z)=\frac{e^{i(a+1)z}}{z+i},\qquad\forall z\in \D.
\]
belongs to $\lambda(p)\setminus\{0\}$.
One can verify this by a simple direct computation and we omit the details.

Thus Theorem~\ref{Small} is proved.
\end{proof}

\begin{remark}
In  this remark, we explain the main idea of the proof of
theorem~\ref{SmallBehaviour}. In general, its feature is demonstrably
application of scheme of Phragm\'en-Lindel\"of theorems generalizing the
maximum modulus principle. Thus, for a given function belonging to
$\lambda(p)$, some auxilliary holomorphic function is constructed and
explored with the  help of the usual maximum modulus principle.
\end{remark}

\begin{proof}[Proof of Theorem~\ref{SmallBehaviour}]
Let $f\in\lambda(p)\setminus\{0\}$ be such that
\[
\liminf_{t\to\infty}\frac{\ln Mf(t)}{t}<\infty.
\]
In particular, from the inclusion $\lambda(p)\subset\Lambda(p)$, it follows
$f\in\Lambda(p)\setminus\{0\}$. Then by Lemma~\ref{Convexity}, there exist
two numbers  $a_0$, $b_0$ such that
\[
\ln Mf(t)\ge a_0t+b_0,\qquad\forall t>0.
\]
Hence $\liminf\limits_{t\to\infty}\frac{\ln Mf(t)}{t}\ge a_0$. We set
\begin{equation}\label{Def_a}
a=\liminf_{t\to\infty}\frac{\ln Mf(t)}{t}.
\end{equation}
So, $a_0\le a<\infty$. In particular,
\[
F(z)=e^{iaz}f(z),\qquad\forall z\in \D,
\]
is a well-defined holomorphic function in  the upper half plane $\D$. We
prove that the function $F$ has the  following properties:
\begin{enumerate}
\item $\ln MF(t)=\ln Mf(t)-at$, $\forall t>0$ --- which  is obvious and we
omit the details;
\item $\ln MF$ is  convex in $(0,\infty)$ --- indeed, $\ln Mf$ is convex in
$(0,\infty)$ by Lemma~\ref{Convexity} and hence, according to the previous
item, $\ln MF$ is convex, too;
\item $\liminf\limits_{t\to\infty}\frac{\ln MF(t)}{t}=0$ --- this  follows
immediately from  both the  equation~(\ref{Def_a})  and the first item;
\item $\ln MF$ is a decreasing function in $(0,\infty)$ --- we deduce this
from both the second and  the third items: for arbitrary $t_1$ and  $t_2$,
$0<t_1<t_2<\infty$, and for $t>t_2$ by the second  item  it follows
\[
\ln MF(t_2)\le\frac{t-t_2}{t-t_1}\ln MF(t_1)+\frac{t_2-t_1}{t-t_1}\ln MF(t)
\]
and we obtain by the third item $\ln MF(t_2)\le\ln MF(t_1)$, i.e. $\ln MF$
decreases in $(0,\infty)$;
\item $\ln MF$ is bounded from above on $[\frac12,\infty)$ --- it follows
immediately from the previous item that $\ln MF(t)\le\ln MF(\frac12)$,
$\forall t\ge\frac12$. Further from $f\in\lambda(p)$ it follows that $\ln
Mf(\frac12)$ is finite and hence by the first item $\ln MF(\frac12)$ is
finite, too, which yields the boundedness of $\ln MF$ on $[\frac12,\infty)$;
\item $F(x+iy)$ tends uniformly  to $0$ as $\abs{x}\to\infty$, with respect
to $y$, $y\in[\frac1c,c]$, where $c\ge 1$is any. We obtain this by the
corresponding property of the function $f\in\lambda(p)$, namely according to
the Definition~\ref{Def}, $f(x+iy)$ has the same property.
\end{enumerate}

We set $A=\sup\limits_{z,\im z\ge1}\abs{F(z)}$. According to the fifth item
$A<\infty$. Moreover, $A\ne 0$ because $F$ is a non zero element of
$\Lambda(p)$.

Note, according to the first item, the assertion  of
Theorem~\ref{SmallBehaviour} can be stated in terms of the function $F$ as
follows
\begin{equation}\label{NewForm}
\lim_{t\to \infty}MF(t)=0.
\end{equation}
In order to prove the  equation~(\ref{NewForm}), we have to obtain that for
every positive number $\varepsilon$ there exists a  number
$y_{\varepsilon}>0$ such that $MF(t)\le\varepsilon$, $\forall
t>y_{\varepsilon}$.

Let $\varepsilon$ be an arbitrary number such that $\varepsilon\in(0,A)$.

It follows  by the sixth item of the above listed properties of the
function $F$ that there  exists a number $x_1$, $x_1>0$, such that
\[
\sup_{x,\abs{x}\ge x_{1}}\abs{F(x+i)}\le\frac{\varepsilon}{2}.
\]
Further, we fix $y_{\varepsilon}$, $y_{\varepsilon}>1$, in such a way  that
\[
\sup_{y,y\ge y_{\varepsilon}}\sqrt{\frac{x_1^2+1}{x_1^2+(y+1)^2}}
\le\frac{\varepsilon}{2A}.
\]

We intend to show that $MF(t)\le\varepsilon$, $\forall t>y_{\varepsilon}$.

Let $z_0\in\D$  be such that $\im z_0>y_{\varepsilon}$.

Let $\eta$, $\eta>0$, be any. We set
\[
y_{\eta}=1+\max\left\{\im
z_0,\frac{1}{\eta}\ln\frac{2A}{\varepsilon}\right\}.
\]
Then it follows from the fifth item that there exists $x_{\eta}$, $$
x_{\eta}>\max\{\abs{\re z_0},x_1\}
$$
such that
\[
\sup_{z,\abs{\re z}\ge x_{\eta},\im z\in
[1,y_{\eta}]}\abs{F(z)}\le\frac{\varepsilon}{2}
\]
We construct the rectangle $\Rect$
\[
\Rect=\{z|\,\abs{\re z}\le x_{\eta}, \im z\in[1,y_{\eta}]\},
\]
and define the holomorphic function $G$ in the upper half plane $\D$ by
\[
G(z)=e^{i{\eta}z}\frac{z}{z+iy_{\varepsilon}}F(z),\qquad\forall z\in \D.
\]

By the choise  of $x_{\eta}$, $y_{\eta}$, $y_{\varepsilon}$ it  is easy to
verify that the modulus $\abs{G}$  of  the  function $G$ is less or equal to
$\frac{\varepsilon}{2}$ on the boundary of the rectangle $\Rect$. Indeed:
\begin{itemize}
\item for $z$ such that $\re z\in[-x_{\eta},x_{\eta}]$, $\im z=y_{\eta}$, we
obtain
\[
\abs{G(z)}=e^{-{\eta}y_{\eta}}\abs{\frac{z}{z+iy_{\varepsilon}}}\abs{F(z)}
\le e^{-{\eta}y_{\eta}}A\le\frac{\varepsilon}{2};
\]
\item for $z$ such that $\re z\in[-x_{\eta},x_{\eta}]$, $\im z=1$, we obtain
\[
\abs{G(z)}=e^{-\eta}\abs{\frac{z}{z+iy_{\varepsilon}}}\abs{F(z)}
\le \abs{\frac{z}{z+iy_{\varepsilon}}}\abs{F(z)}
\]
and further there are  two cases
\begin{enumerate}
\item if $\abs{\re z}\le  x_1$ then we have to continue with the following
inequalities in which we use the  notation $x=\re z$
\[
\begin{array}{rl}
\abs{G(z)}&\displaystyle\le\sqrt{\frac{x^2+1}{x^2+(y_{\varepsilon}+1)^2}
\,}\, \abs{F(z)}\\[6pt]
&\displaystyle\le\sqrt{\frac{x_1^2+1}{x_1^2+(y_{\varepsilon}+1)^2} \,}\,
\abs{F(z)}\\[6pt]
&\displaystyle\le\frac{\varepsilon}{2A}A=\frac{\varepsilon}{2},
\end{array}
\]
\item if $\abs{\re z}\ge x_1$ then it is obvious that
$\abs{G(z)}\le\frac{\varepsilon}{2}$,
\end{enumerate}
i.e. in both cases $\abs{G(z)}\le\frac{\varepsilon}{2}$;
\item for $z$ such that $\im z\in[1,y_{\eta}]$, $\abs{\re z}=x_{\eta}$ we
obtain
\[
\abs{G(z)}=e^{-{\eta}\im z}\abs{\frac{z}{z+iy_{\varepsilon}}}\abs{F(z)}\le
\abs{F(z)}\le\frac{\varepsilon}{2}.
\]
\end{itemize}
Thus, the modulus $\abs{G}$  of  the  function $G$ is less or equal to
$\frac{\varepsilon}{2}$ on the boundary of the rectangle $\Rect$.

Then, by the maximum modulus principle, it follows
\[
\sup_{z,z\in R}\abs{G(z)}\le\frac{\varepsilon}{2},
\]
and hence, in particular, $\abs{G(z_0)}\le\frac{\varepsilon}{2}$.

Therefore,
\[
\begin{array}{rl}
\abs{F(z_0)}&\displaystyle\le e^{{\eta}{\im
z_0}}\abs{\frac{z_0+iy_{\varepsilon}}{z_0}} \frac{\varepsilon}{2}\\[6pt]
&\displaystyle\le  e^{{\eta}{\im
z_0}}\left(1+\frac{y_{\varepsilon}}{\abs{z_0}}
\right)\frac{\varepsilon}{2}\\[6pt]
&\displaystyle\le e^{{\eta}{\im z_0}}\left(1+\frac{y_{\varepsilon}}{\im z_0}
\right)\frac{\varepsilon}{2}\\[6pt]
&\le {\varepsilon}e^{{\eta}{\im z_0}}.
\end{array}
\]
Then $\abs{F(z_0)}\le\varepsilon$ because of the choise of $\eta$.

So, $MF(\im z_0)\le\varepsilon$.

Thus the claim $MF(t)\le\varepsilon$, $\forall t>y_{\varepsilon}$, is
proved.

The proof  of Theorem~\ref{SmallBehaviour} is completed.
\end{proof}

\end{document}